\documentclass[a4paper,12pt]{amsart}

\title[Left Connected]{On the Left Connected Subalgebra of the Descent Algebra of a Coxeter Group of Classical Type}

\author{Linus Hellebrandt}
\address{L.H.: Institute of Algebra and Number Theory, University of Stuttgart, Stuttgart, Germany}
\email{linus.hellebrandt@mathematik.uni-stuttgart.de}

\author{G\"otz Pfeiffer}
\address{G.P.: School of Mathematics,
  Statistics and Applied Mathematics, National University of Ireland,
  Galway, University Road, Galway, Ireland}
\email{goetz.pfeiffer@nuigalway.ie}


\makeatletter
\@namedef{subjclassname@2020}{%
  \textup{2020} Mathematics Subject Classification}
\makeatother


\usepackage{tikz}
\usepackage{tikz-cd}
\usepackage[euler-digits]{eulervm}
\usepackage[margin=1in]{geometry}

\newcommand{\Q}{\mathbb{Q}}
\newcommand{\EE}{\mathfrak{E}}
\newcommand{\supp}{\mathrm{supp}}

\newcommand{\stirlingi}{\genfrac{[}{]}{0pt}{}}
\newcommand{\stirlingii}{\genfrac{\{}{\}}{0pt}{}}

\numberwithin{equation}{section}

\newtheorem{num}{Notation}[section]

\newtheorem{theorem}[num]{Theorem}

\newtheorem{lemma}[num]{Lemma}
\newtheorem{prp}[num]{Proposition}
\newtheorem{cor}[num]{Corollary}

\theoremstyle{definition}
\newtheorem{define}[num]{Definition}

\newtheorem{bsp}[num]{Example}

\newtheorem{bem}[num]{Remark}

\begin{document}

\begin{abstract}
  A Coxeter group of classical type $A_n$, $B_n$ or $D_n$ contains a
  chain of subgroups of the same type.  We show that intersections of
  conjugates of these subgroups are again of the same type, and make
  precise in which sense and to what extent this property is exclusive
  to the classical types of Coxeter groups.  As the main tool for the
  proof we use Solomon's descent algebra.  Using Sterling numbers, we
  express certain basis elements of the descent algebra as polynomials
  and derive explicit multiplication formulas for a commutative
  subalgebra of the descent algebra.
\end{abstract}

\maketitle

\section{Introduction}

For subgroups $J$ and $K$ of a finite group $W$, Mackey's Theorem
describes the restriction to $K$ of a character induced from $J$ to
$W$ in terms of intersections of $W$-conjugates of both $J$ and $K$,
which are parametrized by the $(J, K)$-double cosets of $W$.  In such
a context it is natural to ask: which subgroups of $W$ do occur as
intersections of conjugates of $J$ and $K$?

In~\cite{Harishchandra}, for example, this question arises for a
Coxeter group $W$ of type $B_n$, with Coxeter groups $J, K \leq W$ of
types $B_j$, $B_k$, for some $j, k \leq n$.  In this case the
intersections turn out to be Coxeter groups of type $B_l$, $l \leq
n$. Representatives of the conjugacy classes of these subgroups can be
generated by subsets of the Coxeter generators of $W$, and lie in a
maximal chain of such subgroups of $W$.

In this paper we classify all such chains of subgroups of an
irreducible finite Coxeter group $W$, generated by a set $S$ of simple
reflections.  For a finite Coxeter group, the natural place to study
such questions is the descent algebra.

The descent algebra $\Sigma(W)$ is a remarkable subalgebra of the group
algebra $\Q[W]$, introduced by Solomon~\cite{solomon} in~1976.
It has a linear basis consisting of $2^{|S|}$
elements $x_J$, $J \subseteq S$, and its radical
coincides with the kernel of the homomorphism from
$\Sigma(W)$ into the character ring of $W$ which maps $x_J$ to the
permutation character $\pi_J$ of the action of $W$ on the cosets of
the subgroup $W_J$ generated by $J \subseteq S$.

In this setup, the search for maximal chains of (conjugacy classes of)
subgroups which are closed under intersection leads us to consider
subalgebras $\Q[x_M]$ of $\Sigma(W)$, where $M$ is a maximal subset of
$S$.  We have $\dim \Q[x_M] = \#\{\pi_M(w) : w \in W\}$, since the
minimal polynomial of $x_J$, $J \subseteq S$, is
\begin{align*}
  \prod_{a \in \{\pi_J(w) : w \in W\}} (x - a) \in \Q[x],
\end{align*}
by a theorem of Bonnaf\'e--Pfeiffer~\cite{minpol}.  We classify all
cases where $\Q[x_M]$ is spanned as $\Q$-vector space by a ``native''
basis, i.e., a subset of the basis $\{x_J : J \subseteq S\}$.  It
turns out that this is almost exclusively the case when $W$ is of
classical type $A_n$, $B_n$ or $D_n$ and, as in the above example of
the Coxeter group of type $B_n$, the subset $M$ is ``left-connected''
in $S$, in sense that we will make precise shortly.  We also show
that, if $M$ is left-connected in $S$, the native basis elements $x_J$
of $\Q[x_M]$ are integer polynomials in $x_M$.  For each of the three
classical types of Coxeter groups, we provide explicit formulas for
the structure constants of the algebra $\Q[x_M]$ relative to its
native basis.

We proceed as follows.  In Section~\ref{sec:notation}, we set up
notation for a finite Coxeter group $W$ and its descent algebra
$\Sigma(W)$, and we review some useful properties of Coxeter groups of
classical type.
In Section~\ref{sec:Leftconnected}, we define $W_M$ in $W$ as
\emph{maximally left-connected}, if $W_M$ is of type $X_{n-1}$ in $W$ of type
$X_n$, where $X$ is any of the classical types $A$, $B$ or $D$,
subject to some natural constraints on $n$.  The choice of name is
derived from the position of $M$ in the Coxeter diagram of $(W, S)$.
We show that in these maximally left-connected cases the algebra $\Q[x_M]$
has a native basis.
In Section~\ref{sec:native}, we consider the converse statement and argue that in almost all cases where
$M$ is maximal, but not left-connected in $W$, the algebra $\Q[x_M]$
does not have a native basis.  The only
exceptions are dihedral types and that of $B_1 \times A_1$ in $B_3$. We can now state our main result.
\begin{theorem}\label{thm:main}
  Let $W$ be an irreducible finite Coxeter group with Coxeter
  generators $S$.  Let $J \subseteq S$ be a maximal subset.  Then the subalgebra
  $\mathbb{Q}[x_J]$ of the descent algebra admits a native basis if and only
  if
  \begin{enumerate}
  \item [(i)]  $W$ is of classical type $A_n$, $B_n$ or $D_n$, and $J$ is
  maximally left-connected in $S$;
  \item[(ii)] $W$ is of dihedral type $I_2(m)$ with $m \geq 5$; or
  \item[(iii)]$W$ is of type $B_3$ and $W_J$ is of type $B_1 \times A_1$.
  \end{enumerate}
  Furthermore the native basis consists of integer polynomials in $x_J$ exclusively for \textnormal{(i)}.
\end{theorem}

\section{Notation}\label{sec:notation}

Let $(W, S)$ be a finite Coxeter system, and let $\ell(w)$ denote the minimal
length of $w \in W$ as a word in $S$.  For each subset
$J \subseteq S$, let $W_J$ be the subgroup of $W$ generated by $J$.
The set
\begin{align*}
X_J &= \{ w \in W : \ell(ws) > \ell(w) \text{ for all } s \in J \}
\end{align*}
is
a distinguished transversal of the left cosets $w W_J$ in $W$
consisting of representatives of  minimal length in their coset.
For $K \subseteq S$, the set
$X_K^{-1} = \{ w^{-1} : w \in X_K \}$ is a distinguished transversal
of the right cosets $W_K w$ in $W$, and the set
$X_{JK} = X_J^{-1} \cap X_K$ is a distinguished transversal of the double
cosets $W_J w W_K$ in~$W$.

For $J \subseteq S$, define $x_J = \sum_{w \in X_J} w$.  The subspace
$\Sigma(W) = \sum_{J \subseteq S} \Q x_J$ of the group algebra $\Q W$
is a ring with identity $x_S = 1$, and thus a subalgebra of $\Q W$, called the descent algebra
of $W$.  The elements $x_J$, $J \subseteq S$, are linearly
independent, and thanks to Solomon~\cite{solomon},
\begin{align}\label{eq:aJKL}
  x_J x_K = \sum_{d \in X_{JK}} x_{J^d \cap K} = \sum_{L \subseteq S} a_{JKL} x_L\text,
\end{align}
where $a_{JKL} = \# \{d \in X_{JK} \mid J^d \cap K = L\}$, for $J, K, L \subseteq S$.

For $J \subseteq K \subseteq S$, we denote $X_J^{(K)} = X_J \cap W_K$.
Then $X_J^{(K)}$ is a distinguished left transversal of $W_J$ in
$W_K$, and
$X_J =  X_K X_J^{(K)} = \{vw: v \in X_K,\, w \in X_J^{(K)}\}$, by \cite[2.1.5]{gotz}.
Accordingly, we have
\begin{align}\label{eq:trans}
  x_J = x_K x_J^{(K)}\text,
\end{align}
where   $x_J^{(K)} = \sum_{w \in X_J^{(K)}} w$.

\subsection{Leftconnected subgroups}

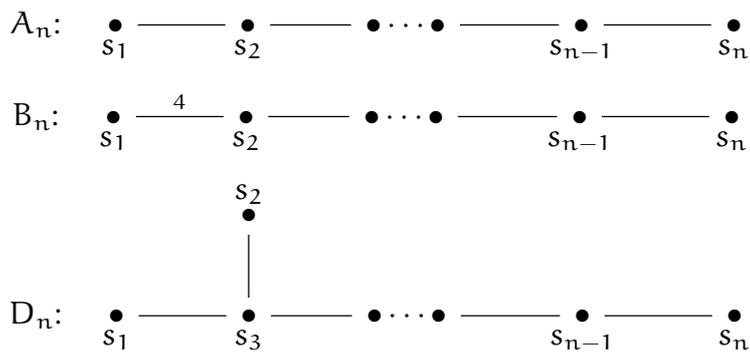
\begin{figure}
\[
\begin{tikzcd}[arrows=-]
A_n{:} &[-2em] \bullet \arrow[r] & \bullet \arrow[r] &  \bullet \cdots \bullet \arrow[r] & \bullet \arrow[r]  & \bullet
\\
[-2.5em] & s_1 & s_2  &    & s_{n-1} & s_n
\end{tikzcd}
\]
\[
\begin{tikzcd}[arrows=-]
B_n{:} &[-2em] \bullet \arrow[r,"4"] & \bullet \arrow[r] &  \bullet \cdots \bullet \arrow[r] &
  \bullet \arrow[r]  & \bullet \\ [-2.5em]
& s_1 & s_2  &    & s_{n-1} & s_n
\end{tikzcd}
\]
\[
\begin{tikzcd}[arrows=-]
&[-2em] &  s_2
\\ [-2.4em]
 &[-2em] &  \bullet \arrow[d]
\\
D_n{:} &[-2em] \bullet \arrow[r] & \bullet \arrow[r] &  \bullet \cdots \bullet \arrow[r] &
  \bullet \arrow[r]  & \bullet \\ [-2.5em]
& s_1 & s_3  &    & s_{n-1} & s_n
\end{tikzcd}
\]
\caption{Coxeter graphs for type $A_n, B_n$ and $D_n$} \label{fig:Coxetergraph}
\end{figure}

In the following section we are mostly concerned with Coxeter groups of type $A_n, B_n$ and $D_n$, whose Coxeter graphs of the respective types are given by Figure \ref{fig:Coxetergraph}.
The types $A_n, B_n$ and $D_n$ are in most parts well understood. For example there exists a well known faithful permutation representation for these types, which in detail is presented in \cite[1.5; 8.1; 8.2]{B_n}:

\begin{bem}
\label{bem:IsomorphismsAndDescents}
(i) For type $A_n$ there exists an isomorphism  between $W$ and the symmetric group $S_{n+1}$ of permutations of the set $\{1,2,\ldots,n+1\}$ which sends the basis element $s_i$ to the involution $(i,i+1)$ interchanging $i$ and $i+1$. This gives rise to an action of $W$ on the set $\{1,2,\ldots,n+1\}$.

(ii) Denote by $S_n^B$ the group of signed permutations on the set $\{-n,-n+1,\ldots,n-1,n\}$, which consists of the permutations $\pi$ satisfying $\pi(-i)=-\pi(i)$ for $0\leq i \leq n$.  For type $B_n$ there exists an isomorphism between $W$ and $S_n^B$ sending $s_i$ to the involution $(-i,-i+1)(i-1,i)$ for $2 \leq i \leq n$, and $s_1$ to $(-1,1)$. This gives rise to an action of $W$ on the set $\{-n,-n+1,\ldots,n-1,n\}$.

(iii) Denote by $S_n^D$ the subgroup of $S_n^B$ consisting of all signed permutations $\pi$ which map an even number of positive integers onto negative integers. For type $D_n$ there exists an isomorphism between $W$ and $S_n^D$ which sends $s_i$ to the involution $(-i,-i+1)(i-1,i)$ for $2 \leq i \leq n$, and $s_1$ to $(-2,1)(-1,2)$.
\end{bem}

Let $W$ be of type $A_n$ ($n \geq 1$), $B_n$ ($n \geq 2$) or $D_n$ ($n \geq 3$), with
$S = \{s_1, s_2, \dots, s_n\}$ as in the diagrams in
Figure~\ref{fig:Coxetergraph}.  For $J = \{s_1, s_2, \dots, s_j\}$ and
$K = \{s_1, s_2, \dots, s_k\}$, $0 \leq j,k \leq n$, we
set $W_j  = W_J$, $X_j = X_J$ and $X_{jk} = X_{JK}$.
Moreover, in $\Q W$, we set $x_j = x_J$ and $x_j^{(k)} = x_J^{(K)}$.
Thus, \eqref{eq:trans} becomes $x_j = x_j^{(k)} x_k$, for $0 \leq j \leq k \leq n$.
For type $D_n$ we exclude the occurrence of $J=\{s_1\}$ or $K=\{s_1\}$ by instead putting $W_1:=W_\emptyset$, $X_1:=X_\emptyset$ and so forth. In Section \ref{subsection:D_n} we furthermore put $x_0:=2x_1$ for type $D_n$, for a simplified notation.

We say that a subgroup $U$ of $W$ is \emph{left connected} in $W$, if
$U = W_j$ for some $0 \leq j \leq n$, where for type $D_n$ we require
$j \geq 1$.  This definition generalizes the notion of left connected
subsets of $S$ for $W$ of type $B_n$, introduced by Gerber, Hi\ss,
Jacon \cite{Harishchandra}, to all classical types of finite Coxeter
groups. We will show in Corollary \ref{cor:combmultAn}, \ref{cor:combmultBn} and \ref{cor:combmultDn} for each of the types $A_n$, $B_n$ and $D_n$ that
\begin{align*}
  x_j x_k = \sum_{l=0}^n a_{jkl} x_l
\end{align*}
for nonnegative integer coefficients $a_{jkl}$, which generalizes \cite[2.2 Lemma]{Harishchandra} for type $A_n$ and $D_n$.

\section{Native bases in Classical types}
\label{sec:Leftconnected}

In this section we define a particular class of subalgebras of $\Sigma(W)$ for any finite Coxeter system. By applying some general properties of the descent algebra to study this algebra for the maximal left connected subgroup in type $A_n, B_n$ and $D_n$, we will see that for this particular case the algebra admits a basis which is a subset of the standard basis of $\Sigma(W)$ and corresponds to a set of falling factorials.

\begin{define}
Let $(W,S)$ be a finite Coxeter system and $\mathfrak{E}\subseteq \Sigma(W)$ be a subalgebra of the descent algebra of $W$. We call a basis of $\mathfrak{E}$ \textit{native} if it is a subset of the standard basis $\{ x_J\mid J\subseteq S\}$ of $\Sigma(W)$. \qed
\end{define}

A trivial example for native bases is the standard basis $\{x_J\mid J\subseteq S\}$ of the descent algebra itself, which is a native basis of $\Sigma(W)$. In the following we will explore non trivial examples of subalgebras which admit native bases.

\begin{define}
\label{define:mathfrakE(W)}
Let $(W,S)$ be a finite Coxeter system and $J\subseteq S$ be a maximal subset of $S$.
Define
\begin{equation*}
\mathfrak{E}(W,x_J):=\mathbb{Q}[x_J],
\end{equation*}
the smallest subalgebra of the descent algebra $\Sigma(W)$, and the group ring $\mathbb{Q}W$, containing the element $x_J$. For type $A_n, B_n$ and $D_n$ put $\mathfrak{E}(W):=\mathfrak{E}(W,x_{n-1})$. \qed
\end{define}

It is clear by its definition that $\mathfrak{E}(W)$ is a commutative algebra isomorphic to a quotient of $\mathbb{Q}[x]$ for all three types via the insertion homomorphism. We now explore the algebras $\mathfrak{E}(W)$ separately for those types. We will see that  $\mathfrak{E}(W)$ admits a native basis and determine its structural constants with respect to this basis. But first we review some properties of falling powers.

\subsection{Falling Factorial Powers}
Following~\cite{concrete}, for $k \geq 0$, the \emph{$k$th falling power} of $x$ is
\begin{align}\label{eq:stirlingi}
  x^{\underline{k}} := \prod_{i=0}^{k-1} (x - i)
  = \sum_{m=0}^k (-1)^{k-m} \stirlingi{k}{m} x^m\text,
\end{align}
where $\stirlingi{k}{m}$ is a Stirling number of the first kind,
the number of ways to arrange $k$ elements in $m$ cycles.  Conversely,
for $k \geq 0$,
\begin{align}\label{eq:stirlingii}
  x^k = \sum_{m=0}^k \stirlingii{k}{m} x^{\underline{m}}\text,
\end{align}
where $\stirlingii{k}{m}$ is a Stirling number of the second kind,
the number of ways to partition an $k$-element set into $m$ nonempty subsets.
Moreover, it is easy to see that
\begin{align}\label{eq:falling-prod}
  x^{\underline{a}}  x^{\underline{b}} = \sum_{k=0}^{\min(a, b)} \binom{a}{k} \binom{b}{k} k!\, x^{\underline{a+b-k}} \text,
\end{align}
where the coefficient $\binom{a}{k} \binom{b}{k} k!$ counts
the number of partial bijections between $k$ elements of a set $A$ with $|A|= a$ and a set $B$ with $|B| = b$.

\subsection{Type $A$}
In the following consider $(W,S)$ a Coxeter system of type $A_n$.
\begin{lemma}
\label{lemma:x_n-1x_n-k}
We have
\begin{align}
\label{eq:A_nx_n-1x_n-k}
x_{n-1}\,x_{n-k}&= k\,x_{n-k}+x_{n-k-1}, \quad \text{ for } 0\leq k \leq n,
\end{align}
where we set $x_{-1}:=x_0$.
So, in particular, $x_{n-1}\, x_0 = n\, x_0 + x_{-1} = (n+1)\, x_0$.
\end{lemma}
\begin{proof}
We show that Equation (\ref{eq:A_nx_n-1x_n-k})
holds for $k<n$, by induction on $k$. For $k=0$ we have $x_{n-k}=x_n=1$, thus $x_{n-1}x_{n-0}=x_{n-1}+0x_{n-0}$. For $k=1$ we have only 2 double cosets of $W_{n-1}$ in $W$ with distinguished coset representatives $1$ and $s_n$. Thus
\begin{equation}
\label{eq:thm:leftconnectedproof}
x_{n-1}^2=x_{n-2}+x_{n-1}
\end{equation} by Equation (\ref{eq:aJKL}), with an exception for $n=1$ where $x_0^2=2x_0$ holds. Thus (\ref{eq:A_nx_n-1x_n-k}) holds for $k\leq 1$, which also proves the case $n\leq 2$.
Suppose $n>2$.
By Equation (\ref{eq:trans}) we have
\begin{equation*}
x_j^{(n)}=x_k^{(n)}x_j^{(k)},
\end{equation*}
for all $0\leq j \leq k\leq n$.
The cases $k\in \{0,1\}$ have already been dealt with, thus let $k\geq 2$ and (\ref{eq:A_nx_n-1x_n-k}) be true for $k-1$. Then
\begin{eqnarray}
\label{eq:inductionleftconnected}
x_{n-1}^{(n)}x_{n-k}^{(n)}&\mathrel{\overset{\makebox[0pt]{\mbox{\normalfont\tiny\sffamily (\ref{eq:trans})}}}{=}} &x_{n-1}^{(n)}x_{n-k+1}^{(n)}x_{n-k}^{(n-k+1)}
\\
\nonumber
&\mathrel{\overset{\makebox[0pt]{\mbox{\normalfont\tiny\sffamily Ind.}}}{=}} &((k-1)x_{n-k+1}^{(n)}+x_{n-k}^{(n)})x_{n-k}^{(n-k+1)}
\\
\nonumber
&\mathrel{\overset{\makebox[0pt]{\mbox{\normalfont\tiny\sffamily (\ref{eq:trans})}}}{=}} & (k-1)x_{n-k}^{(n)}+x_{n-k+1}^{(n)}x_{n-k}^{(n-k+1)}x_{n-k}^{(n-k+1)}
\\
\nonumber
&\mathrel{\overset{\makebox[0pt]{\mbox{\normalfont\tiny\sffamily (\ref{eq:thm:leftconnectedproof})}}}{=}} &(k-1)x_{n-k}^{(n)}+x_{n-k+1}^{(n)}(x_{n-k}^{(n-k+1)}+x_{n-k-1}^{(n-k+1)})
\\
\nonumber
&\mathrel{\overset{\makebox[0pt]{\mbox{\normalfont\tiny\sffamily (\ref{eq:trans})}}}{=}} & kx_{n-k}^{(n)}+x_{n-k-1}^{(n)}.
\end{eqnarray}
Thus by induction (\ref{eq:A_nx_n-1x_n-k}) holds for all $0\leq k <n$, for $n\in \mathbb{N}_+$.
Multiplying $x_{n-1}$ with $x_0$ yields
\begin{equation}
x_{n-1}^{(n)}x_{0}^{(n)}=(n+1)x_0^{(n)},
\end{equation}
since $x_0$ is the sum over all elements in $W$ and $x_{n-1}$ consists of $n+1$ summands, as its number of summands is equal to the index of $W_{n-1}$ in $W$.
\end{proof}

As an immediate consequence of Lemma~\ref{lemma:x_n-1x_n-k},
\begin{align}
\label{eq:identityx_n-kAn}
   x_{n-k} = \prod_{i=0}^{k-1} (x_{n-1} - i) = x_{n-1}^{\underline{k}}\text,\quad
  \text{ for } 0 \leq k \leq n\text,
\end{align}
and
\begin{align*}
(x_{n-1} - n - 1)\, x_0 = x_0 \, (x_{n-1} - n) - x_0 = 0.
\end{align*}
It follows that $x_{n-1}^{\underline{n+1}} = x_{n-1}^{\underline{n}}$
and that $x_{n-1}^{\underline{n+t}} = 0$, for $t > 1$.
This shows the following.

\begin{prp}\label{prp:a}
As $\Q$-algebras  $\EE(W) = \Q[x_{n-1}]$ is isomorphic to $\Q[x]/(x^{\underline{n+1}} - x^{\underline{n}})$. Its $\Q$-bases $\{x_0, x_1, \dots, x_n\}$ and $\{x_{n-1}^0, \dots, x_{n-1}^n\}$ are related by the base change formulas
\begin{align}
  x_{n-k} =  \sum_{m=0}^k (-1)^{k-m} \stirlingi{k}{m} x_{n-1}^m\text,
  \quad \text{ and } \quad
  x_{n-1}^k = \sum_{m=0}^k \stirlingii{k}{m} x_{n-m}\text,
\end{align}
for $0 \leq k \leq n$.\qed
\end{prp}
This confirms that the basis elements $x_j$ of $\EE(W)$ are integral polynomials in $x_{n-1}$,
$0 \leq j \leq n$.
Moreover, the identification with falling powers yields a multiplication formula
in terms of explicit  structure constants of $\mathfrak{E}(W)$ with respect to the basis $\{x_0,x_1,\ldots,x_n\}$.

\begin{cor}
\label{cor:combmultAn}
We have
\begin{eqnarray}
\label{eq:cor:combmult}
x_j\,x_k=\sum_{l=-1}^{\min(j,k)}\binom{n-j}{k-l}\binom{n-k}{j-l}(n-j-k+l)!\,x_l,
\end{eqnarray}
for $0\leq j, k\leq n$.
\end{cor}
\begin{proof}
From \eqref{eq:falling-prod}, we have
\begin{align*}
  x^{\underline{n-j}}\, x^{\underline{n-k}} &= \sum_{s=0}^{\min(n-j,n-k)} \binom{n-j}{n-j-s} \binom{n-k}{n-k-s} s!\, x^{\underline{n-j+n-k-s}} \\
  &= \sum_{l = j+k-n}^{\min(j,k)} \binom{n-j}{k-l} \binom{n-k}{j-l} (n+l-j-k)!\, x^{\underline{n-l}}.
\end{align*}
%where $s = n+l-j-k$, i.e., $l = s-n+j+k$ and $n-j-s = k-l$, and $n-k-s = j - l$.
And, as $x_{n-1}^{\underline{n-l}} = 0$ unless $n-l \leq n + 1$,
it follows that
\begin{align*}
    x_j x_k
  &=
x_{n-1}^{\underline{n-j}}\, x_{n-1}^{\underline{n-k}} =
    \sum_{l = -1}^{\min(j,k)} \binom{n-j}{k-l} \binom{n-k}{j-l} (n+l-j-k)!\, x_l,
\end{align*}
as desired.
\end{proof}

Alternatively we can prove Corollary \ref{cor:combmultAn} using formula (\ref{eq:aJKL})  by counting for each $l\in\{1,2,\ldots, n\}$ the coset representatives $x\in X_{jk}$ satisfying
\begin{equation}
\label{eq:J^xcapK=L}
\{1,2,\ldots,j\}^x\cap \{1,2,\ldots,k\}=\{1,2,\ldots,l\}.
\end{equation}
By \cite[1.5.3]{B_n} we have $\ell(xs_i)>\ell(x)$, if and only if $x(i)<x(i+1)$, which implies \begin{equation*}
 X_j=\{x\in W \mid x(i)<x(i+1), 1 \leq i \leq j\}.
\end{equation*}
Therefore
\begin{align*}
&(i)\hspace{1.5cm} x(1)<x(2)<\cdots<x(k+1),
\\
&(ii) \hspace{1.37cm} x^{-1}(1)<x^{-1}(2)<\cdots<x^{-1}(j+1),
\end{align*}
for each $x\in X_{jk}$.
Thus for a fixed $l$ any coset representative satisfying equation (\ref{eq:J^xcapK=L}) stabilizes the elements $1,2,\ldots, l+1$ point-wise while not stabilizing $l+2$.
The coefficient
\begin{equation*}
\binom{n-j}{n-j-k+l}\binom{n-k}{n-j-k+l}(n-j-k+l)!
\end{equation*}
can then be read as all bijections between the $n-j-k+l$ elements in $\{k+2,\dots, n+1\}$ whose image is not determined by $(ii)$, and the $n-j-k+l$ elements
in $\{j+2,\dots,n+1\}$
in the image whose pre-image is not determined by $(i)$.

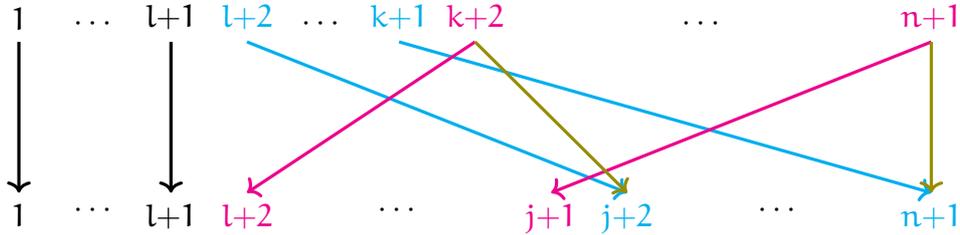
\begin{center}
  \begin{tikzpicture}[scale=2]
    \node[above] (.a) at (0.5,1) {$\cdots$};
    \node[below] (.b) at (0.5,0) {$\cdots$};
    \node[above] (.a) at (2,1) {$\cdots$};
    \node[below] (.b) at (2.5,0) {$\cdots$};
    \node[above] (.a) at (4.5,1) {$\cdots$};
    \node[below] (.b) at (5,0) {$\cdots$};

    \draw[very thick,->] (0,1) node[above] {$1$} -- (0,0) node[below] {$1$};
    \draw[very thick,->] (1,1) node[above] {$l{+}1$} -- (1,0) node[below] {$l{+}1$};
    \draw[cyan,very thick,->] (1.5,1) node[above] {$l{+}2$} -- (4,0) node[below] {$j{+}2$};
    \draw[cyan,very thick,->] (2.5,1) node[above] {$k{+}1$} -- (6,0) node[below] {$n{+}1$};
    \draw[magenta,very thick,->] (3,1) node[above] {$k{+}2$} -- (1.5,0) node[below] {$l{+}2$};
    \draw[magenta,very thick,->] (6,1) node[above] {$n{+}1$} -- (3.5,0) node[below] {$j{+}1$};
    \draw[olive,very thick,->] (6,1) -- (6,0);
    \draw[olive,very thick,->] (3,1) -- (4,0);
  \end{tikzpicture}
\end{center}

\subsection{Type $B$}
In the following consider $(W,S)$ a Coxeter system of type $B_n$.
The set of distinguished double coset representatives of $W_{n-1}$ in $W$ is
\begin{eqnarray*}
X_{n-1,n-1}=\{1,s_n, s_ns_{n-1}\cdots s_2s_1s_2\cdots s_n\},
\end{eqnarray*}
see, e.g., \cite[Example~2.2.5]{gotz}.
Thus $x_{n-1}^2=2x_{n-1}+x_{n-2}$. Analogously to Lemma \ref{lemma:x_n-1x_n-k} induction over $k$ yields:

\begin{lemma}
\label{lemma:x_n-1x-n-kB}
We have
\begin{align}
\label{eq:B_nx_n-1x_n-k}
x_{n-1}x_{n-k}&= 2kx_{n-k}+x_{n-k-1},
\end{align}
for $0\leq k\leq n$, where we put $x_{-1}:=0$. So, in particular, $x_{n-1}x_0=2nx_0+x_{-1}=2nx_0$. \qed
\end{lemma}

As an immediate consequence of Lemma \ref{lemma:x_n-1x-n-kB},
\begin{equation}
\label{eq:x_jaspolynomialinx_n-1B}
x_{n-k}= \prod_{i=0}^{k-1}(x_{n-1}-2i)=2^k \left( \frac{x_{n-1}}{2}\right)^{\underline{k}}, \textnormal{ for $0\leq k \leq n$,}
\end{equation}
and
\begin{equation}
\label{eq:minimalpolynomialB_n}
(x_{n-1}-2n)x_0=0.
\end{equation} This shows the following.

\begin{prp}
\label{prp:isomBn}
As $\Q$-algebras  $\EE(W) = \Q[x_{n-1}]$ is isomorphic to $\Q[x]/(x^{\underline{n+1}})$, where $x_{n-1}$ is sent to the coset of $2x^{\underline{1}}$. Its $\Q$-bases $\{x_0, x_1, \dots, x_n\}$ and $\{x_{n-1}^0, \dots, x_{n-1}^n\}$ are related by the base change formulas
\begin{align}
  x_{n-k} =  \sum_{m=0}^k (-2)^{k-m} \stirlingi{k}{m} x_{n-1}^m\text,
  \quad \text{ and } \quad
  x_{n-1}^k = \sum_{m=0}^k \stirlingii{k}{m} 2^{k-m} x_{n-m},
\end{align}
for $0 \leq k \leq n$.
\end{prp}
\begin{proof}
The minimal polynomial of $x_{n-1}$ arises from Equation (\ref{eq:x_jaspolynomialinx_n-1B}) and (\ref{eq:minimalpolynomialB_n}) as $\mu=\prod_{i=0}^{n}(x-2i)$.
Furthermore the natural isomorphism $\EE(W) \cong \mathbb{Q}[x]/(\mu)$ sends $x_{n-k}$ to the coset of $\prod_{i=0}^{k-1}(x-2i)$. Sending a polynomial $f(x)\in \mathbb{Q}[x]$ to the coset $f(2x)+(x^{\underline{n+1}})$ yields an algebra epimorphism from $\mathbb{Q}[x]$ to $\mathbb{Q}[x]/(x^{\underline{n+1}})$. Since $\mu(2x)=2^{n+1}x^{\underline{n+1}}$ and $x^{\underline{n+1}}$ generate the same ideal of $\mathbb{Q}[x]$ the kernel of this epimorphism is exactly $(\mu)$, thus it induces an isomorphism between $\EE(W)$ and $\mathbb{Q}[x]/(x^{\underline{n+1}})$. Under this isomorphism $x_{n-k}$ is sent to the coset of $2^kx^{\underline{k}}$.
The change of basis formulas of the native basis to the standard basis arise from Equation (\ref{eq:stirlingi}) and (\ref{eq:stirlingii}) under this isomorphism.
\end{proof}

This confirms that the basis elements $x_j$ of $\EE(W)$ are integral polynomials in $x_{n-1}$,
$0 \leq j \leq n$.
Moreover, the identification with multiples of falling powers yields a multiplication formula
in terms of explicit  structure constants of $\mathfrak{E}(W)$ with respect to the basis $\{x_0,x_1,\ldots,x_n\}$. The proof proceeds analogously to Corollary \ref{cor:combmultAn} with the difference that the coset $x^{\underline{n+1}}$ already equals $0$, thus eliminating the second term for $x_0$. This yields:

\begin{cor}
\label{cor:combmultBn}
We have
\begin{equation}
\label{eq:cor:combmult2}
x_j\,x_k=\sum_{l=0}^{\min(j,k)} \binom{n-j}{k-l}\binom{n-k}{j-l}(n-j-k+l)!\,2^{n-j-k+l}x_{l},
\end{equation}
for $0\leq j, k \leq n$.
\qed
\end{cor}

A similar interpretation of the coefficient of $x_i$ as for type $A_n$ exists. Here the action of $W$ on $\{-n,-n+1,\ldots,n\}$ yields $\ell(ws_i)>\ell(w)$, if and only if $w(i-1)<w(i)$, for all $w\in W$ and $1\leq i \leq n$, by \cite[8.1.2]{B_n}. The factor $2^{n-j-k+i}$ arises since for the bijection between the two sets with $n-j-k+i$ elements we can also choose the sign of each element freely.

\subsection{Type $D$}
\label{subsection:D_n}
In the following consider $(W,S)$ a Coxeter system of type $D_n$. The set of distinguished double coset representatives of $W_{n-1}$ in $W$ is
\begin{eqnarray*}
X_{n-1,n-1}=\{1, s_n, s_ns_{n-1}\cdots s_3s_1s_2\cdots s_n\},
\end{eqnarray*}
see, e.g., \cite[Example~2.2.6]{gotz}.
Thus $x_{n-1}^2=2x_{n-1}+x_{n-2}$.  Analogously to Lemma \ref{lemma:x_n-1x_n-k} induction over $k$ yields:
\begin{lemma}
\label{lemma:x_n-1x_n-kD_n}
We have
\begin{align}
\label{eq:D_nx_n-1xn-k}
x_{n-1}x_{n-k}= 2kx_{n-k}+x_{n-k-1}  \text{ for $k\in \{0,1,\ldots,n-1\}$,}
\end{align}
where we set $x_0:=2x_1$. So, in particular, $x_{n-1}x_1=2nx_1$.
\end{lemma}

As an immediate consequence of Lemma \ref{lemma:x_n-1x_n-kD_n},
\begin{equation}
\label{eq:x_jaspolynomialinx_n-1D}
x_{n-k}=
\prod_{i=0}^{k-1}(x_{n-1}-2i)=2^k\left( \frac{x_{n-1}}{2} \right)^{\underline{k}}
, \textnormal{ for $0\leq k \leq n-1$,}
\end{equation}
and
\begin{equation}
\label{eq:minimalpolynomialDn}
(x_{n-1}-2n)x_1=0.
\end{equation} This shows the following.

\begin{prp}
\label{prp:d}
As $\Q$-algebras  $\EE(W) = \Q[x_{n-1}]$ is isomorphic to $\Q[x]/(x^{\underline{n}}-x^{\underline{n-1}})$, where $x_{n-1}$ is sent to the coset of $2x^{\underline{1}}$. Its $\Q$-bases $\{ x_1, \dots, x_n\}$ and $\{x_{n-1}^0, \dots, x_{n-1}^{n-1}\}$ are related by the base change formulas
\begin{align}
  x_{n-k} =  \sum_{m} (-2)^{k-m} \stirlingi{k}{m} x_{n-1}^m\text,
  \quad \text{ and } \quad
  x_{n-1}^k = \sum_m \stirlingii{k}{m} 2^{k-m} x_{n-m},
\end{align}
 for $0 \leq k \leq n-1$.
\end{prp}
\begin{proof}
The minimal polynomial of $x_{n-1}$ arises from Equation (\ref{eq:x_jaspolynomialinx_n-1D}) and (\ref{eq:minimalpolynomialDn}) as $\mu=(x-2n)\prod_{i=0}^{n-2}(x-2i)$. Furthermore the natural isomorphism $\EE(W) \cong \mathbb{Q}[x]/(\mu)$ sends $x_{n-k}$  to the coset of $\prod_{i=0}^{k-1}(x-2i)$. Sending a polynomial $f(x)\in \mathbb{Q}[x]$ to the coset $f(2x)+(x^{\underline{n}}-x^{\underline{n-1}})$ yields an algebra epimorphism from $\mathbb{Q}[x]$ to $\mathbb{Q}[x]/(x^{\underline{n}}-x^{\underline{n-1}})$. Since $\mu(2x)=2^n(x^{\underline{n}}-x^{\underline{n-1}})$ and $x^{\underline{n}}-x^{\underline{n-1}}$ generate the same ideal in $\mathbb{Q}[x]$ the kernel of this epimorphism is exactly $(\mu)$, thus it induces an isomorphism between $\EE(W)$ and $\mathbb{Q}[x]/(x^{\underline{n}}-x^{\underline{n-1}})$. Under this isomorphism $x_{n-k}$ is sent to the coset of $2^kx^{\underline{k}}$. The change of basis formulas of the native basis to the standard basis arise from Equation (\ref{eq:stirlingi}) and (\ref{eq:stirlingii}) under this isomorphism.
\end{proof}

This confirms that the basis elements $x_j$ of $\EE(W)$ are integral polynomials in $x_{n-1}$,
$1 \leq j \leq n$.
Comparing Proposition \ref{prp:d} with Proposition \ref{prp:a} we see that the ideals factoring $\mathbb{Q}[x]$ for type $D_n$ and type $A_{n-1}$ are identical. This gives us a surprising connection between $\EE(W)$ for type $D_n$ and $\EE(W)$ for type $A_{n-1}$.

\begin{cor}
\label{cor:isomorphismDnAn-1}
Let $(W',S')$ be a Coxeter system of type $A_{n-1}$ for $n\geq 3$. To distinguish between elements in $\EE(W)$ and $\EE(W')$ denote by $x_i^D$ the element $x_i\in\EE(W)$ and by $x_i^A$ the element $x_i\in \EE(W')$.
There exists an $\mathbb{Q}$-algebra isomorphism
\begin{equation*}
\varphi: \mathfrak{E}(W)\to\mathfrak{E}(W'),
\end{equation*}
with
\begin{eqnarray*}
\varphi(x_j^D)=
2^{n-j}x_{j-1}^A,
\end{eqnarray*}
for $1\leq j\leq n$. \qed
\end{cor}

We can solve multiplication in $\EE(W)$ by either using the multiplication formula (\ref{eq:falling-prod}) for falling factorials and the isomorphism onto $\mathbb{Q}[x]/(x^{\underline{n}}-x^{\underline{n-1}})$ from Proposition \ref{prp:d}, or using the isomorphism between $\EE(W)$ and $\EE(W')$ in Corollary \ref{cor:isomorphismDnAn-1} and the multiplication formula for type $A_{n-1}$ in Corollary \ref{cor:combmultAn}. Both methods yield the following.

\begin{cor}
\label{cor:combmultDn}
We have
\begin{align}
\label{eq:cor:combmultDn}
x_j\, x_k=
 \sum_{l=0}^{\min(j,k)}\binom{n-j}{k-l}\binom{n-k}{j-l}(n-j-k+l)!\,2^{n-j-k+l}x_{l},
\end{align}
for $1\leq j,k\leq n$, where as a reminder $x_0=2x_1$.
So, in particular, $x_j\,x_1=\frac{n!}{j!}2^{n-j}x_1$.
 \qed
\end{cor}

To conclude our study of $\mathfrak{E}(W)$ we summarize the general results of Proposition \ref{prp:a}, \ref{prp:isomBn} and \ref{prp:d}.

\begin{theorem}\label{thm:classical}
Let $(W,S)$ be a Coxeter system of type $A_n,B_n$ or $D_n$. Then $\EE(W)$ admits a native basis, consisting of all basis elements of $\Sigma(W)$ corresponding to left-connected subgroups of $W$.
Moreover all elements of the native basis are integral polynomials in $x_{n-1}$.
\end{theorem}

We will see in the next section that native bases for $\EE(W,x_J)$ are not exclusive for the leftconnected case. But demanding that the elements of the native basis are integral polynomials in $x_J$ only holds for $\EE(W)$.

\section{No Native Bases}
\label{sec:native}

The following criterion explains why, in most cases of $J$ maximal in $S$,
$\Q[x_J]$ does not have a native basis.

\begin{lemma}\label{la:no-native}
    Let $s \in S$ and set $J:= S \setminus \{s\}$ and $K:= J^s \cap J$.
    \begin{enumerate}
    \item[(i)] If $x_J x_K \neq x_K x_J$ then $\Q[x_J]$ does not have
      a native basis.
    \item [(ii)] If there exist elements $y \in X_{JK}$ and
      $t \in J \setminus K$ such that $t^y \in K$ then
      $x_J x_K\neq x_K x_J$.
\end{enumerate}
\end{lemma}

\begin{proof}
  (i) Assume, for a contradiction, that $\Q[x_J]$ has a native basis.
  This basis contains $x_J$ and $x_K$, as both $J, K \in \supp(x_J^2)$
  (since by \eqref{eq:aJKL}, $x_J^2 = \sum_{x \in X_{JJ}} x_{J^x \cap J}$ and
  clearly $1, s \in X_{JJ}$).  But $x_J x_K \neq x_K x_J$ is absurd
  since $\Q[x_J]$ is commutative.   Therefore
  $\Q[x_J]$ does not have a native basis.

  (ii) Set $x = y^{-1} \in X_{KJ}$.  Then $t \in K^x \cap J$ but
  $t \notin K$.  Hence $K^x \cap J \notin \supp(x_J x_K)$ as
  clearly $L \subseteq K$ for all $L \in \supp(x_J x_K)$.  But
  $K^x \cap J \in \supp(x_K x_J)$ as $x \in X_{KJ}$.  So
  $\supp(x_J x_K) \neq \supp(x_K x_J)$, whence $x_J x_K \neq x_K x_J$.
\end{proof}

We list some explicit cases in Table~\ref{tab:no-native}.
Here, for $J \subseteq K \subseteq S$, we write $d_J^K$
for the longest right coset representative $w_J w_K$ of $W_J$
in $W_K$.  With this notation, in each case, it is straightforward to verify
that the listed element $y$ satisfies the conditions $y \in X_{JK}$ and $t^y \in K$.  Here, and in the remainder of this section,
the labelling of the elements of $S$ follows the conventions used,
e.g., in~\cite[1.3]{gotz}.

In fact,
Table~\ref{tab:no-native} almost suffices to classify all native bases, as follows.

\begin{table}%\setlength\extrarowheight{1.5pt}
  \renewcommand{\arraystretch}{1.2}
  \begin{tabular}{cccccc}
    $W$ & $s$ & $K$ & $t \in J \setminus K$ & $y \in X_{JK}$ & $t^y$ \\
    \hline
    $B_3$ & $s_1$ & $\{s_3\}$ & $s_2$ &
    % $_{(121)(32)} = d_2^{12} d_2^{23}$ or
    $d_{23}^{123}$ & $s_3$ \\
    $H_3$ & $s_1$ & $\{s_3\}$ & $s_2$ &
    % $_{121232121321} =
    $d_{23}^{123}$ & $s_3$ \\
    $H_3$ & $s_3$ & $\{s_1\}$ & $s_2$ &
    % $_{3212132123} =
    $d_{12}^{123}$ & $s_1$ \\
    $A_4$ & $s_2$ & $\{s_4\}$ & $s_1$ &
    % $_{(2132)(43)}  =
    $d_{13}^{123} d_{13}^{134}$ & $s_4$ \\
    $B_4$ & $s_2$ & $\{s_4\}$ & $s_3$ & %$_{(2123212)(43)}  =
    $d_{13}^{123} d_{13}^{134}$ & $s_4$ \\
    $H_4$ & $s_2$ & $\{s_4\}$ & $s_3$ & %$_{(2121321213212)(43)}  =
    $d_{13}^{123} d_{13}^{134}$ & $s_4$ \\
    $H_4$ & $s_4$ & $\{s_1,s_2\}$ & $s_3$ & %$_{(43)(2121)(32)(43)} =
    $d_3^{34} d_1^{12} d_2^{23} d_3^{34}$ & $s_2$ \\
    $D_5$ & $s_1$ & $\{s_2,s_4,s_5\}$ & $s_3$ & %$_{(132431)(543)} =
    $d_{234}^{1234} d_{34}^{345}$ & $s_4$ \\
    $E_6$ & $s_6$ & $\{s_1,s_2,s_3,s_4\}$ & $s_5$ & %$_{(65423456)(1345)} =
    $d_{2345}^{23456} d_{345}^{1345}$ & $s_4$  \\
    $E_7$ & $s_7$ & $\{s_1,s_2,\dots,s_5\}$ & $s_6$ & %$_{(7654234567)(13456)} =
    $d_{23456}^{234567} d_{3456}^{13456}$ & $s_5$  \\
    $E_8$ & $s_8$ & $\{s_1,s_2,\dots,s_6\}$ & $s_7$ & %$_{(876542345678)(134567)} =
    $d_{234567}^{2345678} d_{34567}^{134567}$ & $s_6$  \\
    \hline
  \end{tabular}
  \caption{No Native Basis}\label{tab:no-native}
\end{table}

\begin{prp}\label{prp:native}
  Let $W$ be an irreducible  finite Coxeter group with simple reflections $S$.
  Let $s \in S$ and $J:= S \setminus \{s\}$ be such that $J$ is not left-connected in $S$.  Then the subalgebra
  $\Q[x_J]$ of $\Sigma(W)$ has a native basis if and only if
  \begin{enumerate}
  \item [(i)] $W$ is of type $I_2(m)$ with $m \geq 5$; or
  \item[(ii)] $W$ is of type $B_3$ and $s = s_2$.
  \end{enumerate}
\end{prp}

\begin{proof}
Note that if $L \subseteq S$ is such that $s \in L$ and the
 subgroup $W_L$ of $W$ together with the reflection $s$
occurs as one of the cases in Table~\ref{tab:no-native} then, with the
same notation $J = S \setminus \{s\}$ and $K = J^s \cap J$, it still
follows that $x_J x_K \neq x_K x_J$.  This is because one can choose
the same elements $t \in J \setminus K$ and $y \in X_{JK}$ as in the
table so that $t^y \in K$.  Hence, in these cases, $\Q[x_J]$ does not
have a native basis by Lemma~\ref{la:no-native}.

A careful inspection of the classification of finite irreducible
Coxeter groups shows that, apart from rank $2$ and the left-connected
maximal subsets $J \subseteq S$ this leaves only three cases to
consider: $W$ of type $B_3$ with $s = s_2$, $W$ of type $H_3$ with
$s = s_2$ and $W$ of type $F_4$ with $s = s_1$.  Case (i), where $W$
is of rank~$2$, is described in detail in Example~\ref{ex:rank2}.
Case (ii), where $W$ is of type $B_3$ and $s = s_2$ is described in
detail in Example~\ref{ex:b3}.

In both remaining cases, it can be explicitly verified that
$x_J x_K \neq x_K x_J$ for some $K \in \supp(x_J^2)$:
\begin{itemize}
\item $W = H_3$, $s = s_2$: $J = \{s_1,s_3\}$, $x = d_1^{12} d_2^{23}$, $J^x \cap J = K = \{s_3\}$.  Then $K^y \cap J = \{s_1\} \not\subseteq K$ for $y = s_2s_1s_3s_2s_1s_2$.
\item $W = F_4$, $s = s_1$: $J = \{s_2,s_3,s_4\}$, $x = d_{23}^{123}$, $J^x \cap J = K = \{s_2,s_3\}$. Then $K^y \cap J = \{s_4\} \not\subseteq K$ for $y = s_1s_2s_3s_2s_4s_3s_2s_1$.
\end{itemize}
This concludes the proof of Proposition~\ref{prp:native}.
\end{proof}

\begin{bsp}\label{ex:rank2}
If $W$ is of type $I_2(m)$, $m = 2k + l$ (for $k \geq 1$, $l = 1,2$), with $s = s_1$ and $J = \{s_2\}$,
we have
\begin{align*}
  x_J^0 &= x_S & x_S &= x_J^0\\
  x_J^1 &= x_J & x_J&= x_J^1 \\
  x_J^2 &= lx_J + k x_{\emptyset} & x_{\emptyset} &= -\tfrac{l}k x_J^1 + \tfrac1k x_J^2
\end{align*}
Thus, $\Q[x_J]$ has a native basis but, unless $k = 1$, not all of its elements are
integer polynomials in $x_J$.  If $k = 1$ then $W$ is of
type $A_2$ ($l = 1$) or $B_2$ ($l = 2$) and $J$ is left-connected in $S$.
\end{bsp}

\begin{bsp}\label{ex:b3}
If $W$ is of type $B_3$, with $s = s_2$, $J = \{s_1, s_3\}$ and $K = \{s_1\}$, we have
\begin{align*}
  x_J^0 &= x_S &                      x_S &= x_J^0 \\
  x_J^1 &= x_J &                      x_J &= x_J^1 \\
  x_J^2 &= 2 x_J + x_K + 2 x_{\emptyset} & x_K &= -\tfrac{14}{5} x_J^1 + \tfrac{8}{5} x_J^2 - \tfrac{1}{10} x_J^3 \\
  x_J^3 & = 4x_J + 6 x_K + 32 x_{\emptyset} & x_{\emptyset} &= \tfrac{2}{5} x_J^1 - \tfrac{3}{10} x_J^2 + \tfrac{1}{20} x_J^3
\end{align*}
Thus, $\Q[x_J]$ has a native basis, but not all of its elements are integer polynomials in $x_J$.
\end{bsp}

Our Main Theorem~\ref{thm:main} now follows from Theorem~\ref{thm:classical} and Proposition~\ref{prp:native}.

%%%%%%%%%%%%%%%%%%%%%%%%%%%%%%%%%%%%%%%%%%%%%%%%%%%%%%%%%%%%%%%%%%%%%%%%%%%%%
\bibliography{paper}

\providecommand{\bysame}{\leavevmode\hbox to3em{\hrulefill}\thinspace}
\providecommand{\MR}{\relax\ifhmode\unskip\space\fi MR }
\providecommand{\MRhref}[2]{%
  \href{http://www.ams.org/mathscinet-getitem?mr=#1}{#2}
}
\providecommand{\href}[2]{#2}
\begin{thebibliography}{1}

\bibitem{B_n}
Anders Bj\"{o}rner and Francesco Brenti, \emph{Combinatorics of {C}oxeter
  groups}, Graduate Texts in Mathematics, vol. 231, Springer, New York, 2005.
  \MR{2133266}

\bibitem{minpol}
C.~Bonnaf\'{e} and G.~Pfeiffer, \emph{Around {S}olomon's descent algebras},
  Algebr. Represent. Theory \textbf{11} (2008), no.~6, 577--602. \MR{2453230}

\bibitem{gotz}
Meinolf Geck and G\"{o}tz Pfeiffer, \emph{Characters of finite {C}oxeter groups
  and {I}wahori-{H}ecke algebras}, London Mathematical Society Monographs. New
  Series, vol.~21, The Clarendon Press, Oxford University Press, New York,
  2000. \MR{1778802}

\bibitem{Harishchandra}
Thomas Gerber, Gerhard Hiss, and Nicolas Jacon, \emph{Harish-{C}handra series
  in finite unitary groups and crystal graphs}, Int. Math. Res. Not. IMRN
  (2015), no.~22, 12206--12250. \MR{3456719}

\bibitem{concrete}
Ronald~L. Graham, Donald~E. Knuth, and Oren Patashnik, \emph{Concrete
  mathematics}, second ed., Addison-Wesley Publishing Company, Reading, MA,
  1994, A foundation for computer science. \MR{1397498}

\bibitem{solomon}
Louis Solomon, \emph{A {M}ackey formula in the group ring of a {C}oxeter
  group}, J. Algebra \textbf{41} (1976), no.~2, 255--264. \MR{444756}

\end{thebibliography}
\bibliographystyle{amsplain}

\end{document}